\documentclass[10pt]{amsart}

\usepackage{amssymb,amsmath,amsthm,amsfonts}

\input{prepictex}
\input{pictex}
\input{postpictex}

\theoremstyle{plain}
\newtheorem{Theorem}{Theorem}
\newtheorem{Proposition}[Theorem]{Proposition}
\newtheorem{Lemma}[Theorem]{Lemma}
\newtheorem{Corollary}[Theorem]{Corollary}

\theoremstyle{definition}

\begin{document}
\title{Counting $(1,\beta)$-BM relations \\ and \\ classifying $(2,2)$-BM groups}
\author{Diego Rattaggi}
\email{rattaggi@hotmail.com}
\date{July 5, 2006}
\begin{abstract}
In the first part, we prove that the number of $(1,\beta)$-BM 
relations is $3 \cdot 5 \cdot \ldots \cdot (2\beta + 1)$, which
was conjectured by Jason Kimberley.
In the second part, we construct two isomorphisms between certain $(2,2)$-BM groups.
This completes the classification of $(2,2)$-BM groups initiated in \cite{KR}.
\end{abstract}
\maketitle

\section{Introduction} \label{Intro}
Let $\mathcal{T}_{r}$ be the $r$-regular tree and
$\mathrm{Aut}(\mathcal{T}_{r})$ its group of automorphisms.
If $\alpha, \beta \in \mathbb{N}$,
an \emph{$(\alpha,\beta)$-BM group} is a torsion-free subgroup of
$\mathrm{Aut}(\mathcal{T}_{2\alpha}) \times \mathrm{Aut}(\mathcal{T}_{2\beta})$
acting freely and transitively on the vertex set of the affine 
building $\mathcal{T}_{2\alpha} \times \mathcal{T}_{2\beta}$.

The class of $(\alpha,\beta)$-BM groups includes for example
$F_{\alpha} \times F_{\beta}$ (the direct product of free groups
of rank $\alpha$ and $\beta$), but also more complicated groups,
like groups containing a finitely presented, torsion-free, simple
subgroup of finite index, if $\alpha$ and $\beta$ are large enough, see \cite[Theorem~6.4]{BMII}.
The first (and only known) examples of
finitely presented, torsion-free, simple groups
have been found in this way.
See also \cite[Section~II.5]{Wise} 
for a non-residually finite $(4,3)$-BM group, 
\cite[Example~2.3]{Rattaggi2} for a $(3,3)$-BM group
having no non-trivial normal subgroups of infinite index,
and \cite[Example~3.4]{Rattaggi2} for a $(6,4)$-BM group
having a subgroup of index $4$ which is finitely presented,
torsion-free, and simple.

\begin{figure}[htbp]
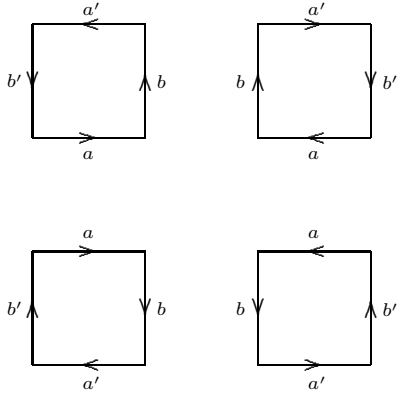
\label{fig1}
\hfil
\centerline{
\beginpicture
\setcoordinatesystem units <0.75cm, 0.755cm>  
\setplotarea  x from 0 to 6,  y from 0 to 6
\putrule from 0 4 to 2 4
\putrule from 0 6 to 2 6
\putrule from 2 4 to 2 6
\putrule from 0 4 to 0 6
\arrow <6pt> [.3,.67] from   1.09  4   to  1.11 4 
\put {$_{a}$}      at   1  3.7
\arrow <6pt> [.3,.67] from   2  5.09   to  2 5.11
\put {$_{b}$}      at   2.3  5
\arrow <6pt> [.3,.67] from   0.91  6   to  0.89 6 
\put {$_{\phantom{'}a'}$}      at   1  6.3
\arrow <6pt> [.3,.67] from   0  4.91   to  0 4.89
\put {$_{b'}$}      at   -0.3  5
\putrule from 4 4 to 6 4
\putrule from 4 6 to 6 6
\putrule from 6 4 to 6 6
\putrule from 4 4 to 4 6
\arrow <6pt> [.3,.67] from   4.91  4   to  4.89 4 
\put {$_{a}$}      at   5  3.7
\arrow <6pt> [.3,.67] from   6  4.91   to  6 4.89
\put {$_{\phantom{'}b'}$}      at   6.3  5
\arrow <6pt> [.3,.67] from   5.09  6   to  5.11 6 
\put {$_{\phantom{'}a'}$}      at   5  6.3
\arrow <6pt> [.3,.67] from   4  5.09   to  4 5.11
\put {$_{b}$}      at   3.7  5
\putrule from 0 0 to 2 0
\putrule from 0 2 to 2 2
\putrule from 2 0 to 2 2
\putrule from 0 0 to 0 2
\arrow <6pt> [.3,.67] from   0.91  0   to  0.89 0 
\put {$_{\phantom{'}a'}$}      at   1  -0.3
\arrow <6pt> [.3,.67] from   2  0.91   to  2 0.89
\put {$_{b}$}      at   2.3  1
\arrow <6pt> [.3,.67] from   1.09  2   to  1.11 2 
\put {$_{a}$}      at   1  2.3
\arrow <6pt> [.3,.67] from   0  1.09   to  0 1.11
\put {$_{b'}$}      at   -0.3  1
\putrule from 4 0 to 6 0
\putrule from 4 2 to 6 2
\putrule from 6 0 to 6 2
\putrule from 4 0 to 4 2
\arrow <6pt> [.3,.67] from   5.09  0   to  5.11 0 
\put {$_{\phantom{'}a'}$}      at   5  -0.3
\arrow <6pt> [.3,.67] from   6  1.09   to  6 1.11
\put {$_{\phantom{'}b'}$}      at   6.3  1
\arrow <6pt> [.3,.67] from   4.91  2   to  4.89 2 
\put {$_{a}$}      at   5  2.3
\arrow <6pt> [.3,.67] from   4  0.91   to  4 0.89
\put {$_{b}$}      at   3.7  1
\endpicture
}
\hfil
\caption{The geometric square $[aba'b']$, represented by each of these
four squares.}
\end{figure}

An equivalent definition for an $(\alpha,\beta)$-BM group is the following 
(the equivalence is shown in \cite[Theorem~3.4]{KR}):
Let $A_{\alpha} = \{ a_1, \ldots, a_{\alpha} \}$,
$B_{\beta} = \{ b_1, \ldots, b_{\beta} \}$, $a, a' \in A_{\alpha}^{\pm 1}$, 
and $b, b' \in B_{\beta}^{\pm 1}$.
We think of the elements in $A_{\alpha}^{\pm 1}$ as oriented horizontal edges
and the elements in $B_{\beta}^{\pm 1}$ as oriented vertical edges.
A \emph{geometric square} $[aba'b']$ is 
a set (consisting of a usual oriented square $aba'b'$ and reflections along
its edges)
\[
[aba'b'] := \{a b a' b', \, a'b'ab, \, a^{-1}b'^{-1}a'^{-1}b^{-1}, \,
a'^{-1}b^{-1}a^{-1}b'^{-1} \}.
\]
See Figure~\ref{fig1} for an illustration of the geometric square $[aba'b']$.

It is easy to check that
\[
[a b a' b'] = [a'b'ab] = [a^{-1}b'^{-1}a'^{-1}b^{-1}] =
[a'^{-1}b^{-1}a^{-1}b'^{-1}].
\]
Let $GS_{\alpha, \beta}$ be the set of all such geometric squares.
\[
GS_{\alpha, \beta} := \{[aba'b'] : a, a' \in A_{\alpha}^{\pm 1}, \, b, b' \in B_{\beta}^{\pm 1} \}.
\]
Given a subset $S \subseteq GS_{\alpha, \beta}$, the \emph{link} $Lk(S)$ is an 
undirected graph with vertex set 
$A_{\alpha}^{\pm 1} \sqcup B_{\beta}^{\pm 1}$
and edges $\{a^{-1}, b \}$, $\{a', b^{-1} \}$, $\{a'^{-1}, b' \}$,
$\{a, b'^{-1} \}$ for each geometric square $[aba'b'] \in S$.
These edges in the link correspond to the four corners
in $aba'b'$.
An \emph{$(\alpha,\beta)$-BM relation} is a
set $R$ consisting of exactly $\alpha \beta$ geometric squares in $GS_{\alpha, \beta}$
such that $Lk(R)$ is the complete bipartite graph
$K_{2\alpha, 2\beta}$ (where the bipartite structure is 
induced by the decomposition 
$A_{\alpha}^{\pm 1} \sqcup B_{\beta}^{\pm 1}$).
This \emph{link condition} for $R$ means that for any given
$a \in A_{\alpha}^{\pm 1}$, $b \in B_{\beta}^{\pm 1}$,
there are unique $a' \in A_{\alpha}^{\pm 1}$, 
$b' \in B_{\beta}^{\pm 1}$ such that $[aba'b'] \in R$.
It also excludes the existence of
geometric squares of the form $[abab]$ in an $(\alpha,\beta)$-BM
relation by a simple counting argument 
($K_{2\alpha, 2\beta}$ has $2\alpha + 2 \beta$ vertices
and $2 \alpha \cdot 2 \beta = 4 \alpha \beta$ edges, 
so each of the $\alpha \beta$ geometric squares
in $R$ has to contribute four distinct edges,
but $[abab]$ only contributes the two edges $\{ a^{-1}, b \}$ and $\{ a, b^{-1} \}$).
We denote by $R_{\alpha, \beta}$ the set of
$(\alpha,\beta)$-BM relations.
Any group $\Gamma$ with a finite presentation 
$\langle A_{\alpha} \cup B_{\beta} \mid R \rangle$, where $R \in R_{\alpha, \beta}$,
is called an \emph{$(\alpha,\beta)$-BM group}. 
Note that any of the four squares representing a geometric square
induces the same relation in $\Gamma$, and that therefore
any $(\alpha,\beta)$-BM group has a presentation with
$\alpha + \beta$ generators and $\alpha \beta$ relations
of the form $aba'b'$.

The cardinality of $R_{\alpha, \beta}$ 
(i.e.\ the number of $(\alpha,\beta)$-BM relations)
has been computed for a finite number of small pairs $(\alpha, \beta)$
in \cite[Table~B.3]{Rattaggi}
and independently with a different method in \cite[Table~4]{JSK},
see Table~\ref{table1}.
\begin{table}[ht] \label{table1}
\[
\begin{tabular}{| c | c | r l |}
\hline
$\alpha$ & $\beta$ &  $|R_{\alpha, \beta}|$ & \\ \hline
$1$ & $1$ & $3$ & \\
$1$ & $2$ & $15$ & $ = 3 \cdot 5$ \\
$1$ & $3$ & $105$ & $ = 3 \cdot 5 \cdot 7$ \\
$1$ & $4$ & $945$ & $ = 3 \cdot 5 \cdot 7 \cdot 9$ \\
$1$ & $5$ & $10395$ & $ = 3 \cdot 5 \cdot 7 \cdot 9 \cdot 11$ \\
$1$ & $6$ & $135135$ & $ = 3 \cdot 5 \cdot 7 \cdot 9 \cdot 11 \cdot 13$ \\
$1$ & $7$ & $2027025$ & $ = 3 \cdot 5 \cdot 7 \cdot 9 \cdot 11 \cdot 13 \cdot 15$ \\
$1$ & $8$ & $34459425$ & $ = 3 \cdot 5 \cdot 7 \cdot 9 \cdot 11 \cdot 13 \cdot 15 \cdot 17$ \\
$1$ & $9$ & $654729075$ & $ = 3 \cdot 5 \cdot 7 \cdot 9 \cdot 11 \cdot 13 \cdot 15 \cdot 17 \cdot 19$ \\
$2$ & $2$ & $541$ & prime\\
$2$ & $3$ & $35235$ & $ = 3^5 \cdot 5 \cdot 29$\\
$2$ & $4$ & $3690009$ & $ = 3^3 \cdot 19 \cdot 7193$\\
$2$ & $5$ & $570847095$ & $ = 3^6 \cdot 5 \cdot 7 \cdot 13 \cdot 1721$\\
$3$ & $3$ & $27712191$ & $ = 3 \cdot 13 \cdot 710569$\\
\hline
\end{tabular}
\]
\caption{Number of $(\alpha,\beta)$-BM relations, $\alpha \leq \beta$.} \label{TableR}
\end{table}

In the smallest case, we have $|R_{1,1}| = 3$, since
\[
R_{1,1} = \big\{ \{ [a_1 b_1 a_1^{-1} b_1^{-1}] \}, 
\{ [a_1 b_1 a_1 b_1^{-1}] \}, 
\{ [a_1 b_1 a_1^{-1} b_1] \} \big\},
\]
using the observation that 
\[
\{ [a_1 b_1 a_1 b_1] \} = \{ [a_1^{-1} b_1^{-1} a_1^{-1} b_1^{-1}] \} \notin R_{1,1}
\]
and 
\[
\{ [a_1 b_1^{-1} a_1 b_1^{-1}] \} = \{ [a_1^{-1} b_1 a_1^{-1} b_1] \} \notin R_{1,1}.
\]
In general, let's say if $\alpha \beta > 10$,
the value $|R_{\alpha, \beta}|$ is not known, but
Kimberley has conjectured in \cite[Conjecture~193]{JSK} that 
$|R_{1, \beta}| = 3 \cdot 5 \cdot \ldots \cdot (2\beta + 1)$ 
for all $\beta \in \mathbb{N}$.
We will prove this conjecture in Section~\ref{SectionCount}.
Observe that $|R_{\alpha, \beta}| = |R_{\beta, \alpha}|$
and therefore $|R_{\alpha, 1}| = 3 \cdot 5 \cdot \ldots \cdot (2\alpha + 1)$ 
for all $\alpha \in \mathbb{N}$.

Each element $R \in R_{\alpha, \beta}$ defines the
$(\alpha,\beta)$-BM group 
$\langle A_{\alpha} \cup B_{\beta} \mid R \rangle$.
Of course, it is possible that distinct $(\alpha,\beta)$-BM relations
define isomorphic $(\alpha,\beta)$-BM group, for example 
(taking $\alpha = \beta = 1$) 
\[
\langle a_1, b_1 \mid a_1 b_1 a_1 b_1^{-1} \rangle \cong 
\langle a_1, b_1 \mid a_1 b_1 a_1^{-1} b_1 \rangle,
\]
whereas $\{ [a_1 b_1 a_1 b_1^{-1}] \} \ne \{[a_1 b_1 a_1^{-1} b_1] \}$.
The classification of 
$(\alpha,\beta)$-BM groups up to isomorphism seems
to be a hard problem in general 
(even if the set $R_{\alpha, \beta}$ is known). 
It has been done
by Kimberley in \cite[Chapter~5]{JSK} for $(1,\beta)$-BM groups,
if $\beta \in \{ 1, \ldots, 5 \}$.
Moreover, Kimberley and Robertson have proved that there are at least $41$
and at most $43$ $(2,2)$-BM groups up to isomorphism,
see \cite[Section~7]{KR} and \cite[Chapter~5]{JSK}.
Starting from a reservoir of $|R_{2,2}| = 541$ $(2,2)$-BM relations,
the lower bound was achieved by computing the
abelianizations of the corresponding $(2,2)$-BM groups, 
and the abelianizations of subgroups of low index.
The upper bound comes from constructing isomorphisms via
generator permutations and Tietze transformations.
It remained the open question whether the group $\Gamma_{4}$ is isomorphic to $\Gamma_{30}$
and whether $\Gamma_{5}$ is isomorphic to $\Gamma_{10}$
(these four $(2,2)$-BM groups will be defined in Section~\ref{SectionClass}).
We will give a positive answer by constructing explicit isomorphisms,
such that there are in fact exactly $41$ $(2,2)$-BM groups up to isomorphism.
If $\alpha, \beta \geq 2$, no other complete classification of 
$(\alpha,\beta)$-BM groups is known so far.

\section{Counting $(1,\beta)$-BM relations} \label{SectionCount}
In this section, we will define a map $\psi_{\beta}$ which
associates to any $(1,\beta)$-BM relation $R \in R_{1,\beta}$
a set $\psi_{\beta}(R) = \psi_{\beta}^{(1)}(R) \cup \psi_{\beta}^{(2)}(R)$ 
consisting of $3 + 2\beta$ distinct
$(1,\beta + 1)$-BM relations 
(see Lemma~\ref{Lemma2} and Lemma~\ref{Lemma3}). 
These $3 + 2\beta$ elements
are either obtained by adding to $R$ a single new
geometric square, or by first removing from $R$
one of the $\beta$ geometric squares and then adding two suitably chosen
new geometric squares. Distinct elements $R$, $T$ in $R_{1,\beta}$
will produce disjoint sets $\psi_{\beta}(R)$, $\psi_{\beta}(T)$ (see Lemma~\ref{Lemma4}).
Moreover, any $(1,\beta + 1)$-BM relation can be obtained by $\psi_{\beta}$
(see Lemma~\ref{Lemma5}).
This allows us to compute inductively the exact number of
$(1,\beta)$-BM relations for any $\beta \in \mathbb{N}$,
and therefore to prove Kimberley's conjecture.

Let
$R = \{r_1, \ldots , r_{\beta}\} \in R_{1,\beta}$,
i.e.\ $r_1, \ldots , r_{\beta}$ are $\beta$ geometric squares
in $GS_{1,\beta}$ satisfying the link condition 
$Lk(\{r_1, \ldots , r_{\beta}\}) = K_{2,2 \beta}$.
We first define
\begin{align}
\psi_{\beta}^{(1)}(R) := 
\big\{ &\{r_1, \ldots, r_{\beta},[a_1 b_{\beta + 1} a_1^{-1} b_{\beta + 1}^{-1}]\},  \notag \\
&\{r_1, \ldots, r_{\beta},[a_1 b_{\beta + 1} a_1 b_{\beta + 1}^{-1}]\}, \notag \\
&\{r_1, \ldots, r_{\beta},[a_1 b_{\beta + 1} a_1^{-1} b_{\beta + 1}]\}
\big\}, \notag
\end{align}
a set consisting of three distinct
$(1,\beta + 1)$-BM relations.

If $[aba'b'] \in GS_{1,\beta}$, we define
\[
\phi_{\beta}([aba'b']) := \big\{ \{ [a b_{\beta + 1} a' b'], [a b a'  b_{\beta + 1}^{-1}] \}, 
\{[a b_{\beta + 1}^{-1} a' b'], [a b a' b_{\beta + 1}] \} \big\}. 
\]

\begin{figure}[htbp]
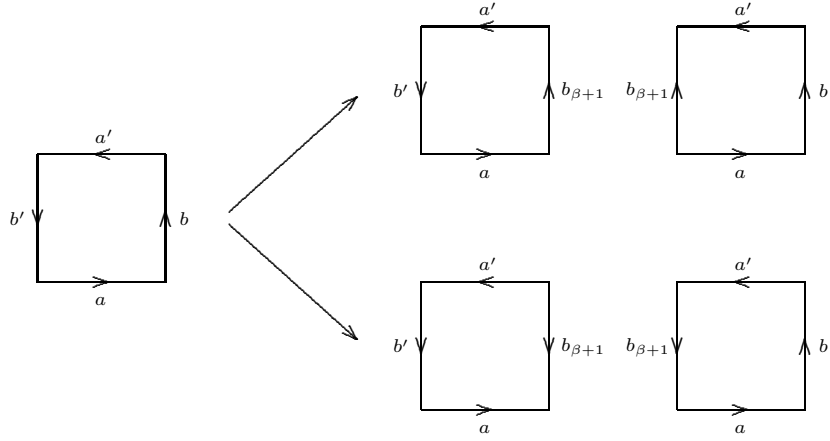
\label{fig2}
\hfil
\centerline{
\beginpicture
\setcoordinatesystem units <0.85cm, 0.85cm>  
\setplotarea  x from 0 to 12,  y from 0 to 6
\putrule from 6 4 to 8 4
\putrule from 8 4 to 8 6
\putrule from 6 6 to 8 6
\putrule from 6 4 to 6 6
\arrow <6pt> [.3,.67] from   7.09  4   to  7.11 4 
\put {$_{a}$}      at   7  3.7
\arrow <6pt> [.3,.67] from   8  5.09   to  8 5.11
\put {$_{b_{\beta + 1}\phantom{'}}$}      at   8.6  5
\arrow <6pt> [.3,.67] from   6.91  6   to  6.89 6 
\put {$_{\phantom{'}a'}$}      at   7  6.3
\arrow <6pt> [.3,.67] from   6  4.91   to  6 4.89
\put {$_{b'}$}      at   5.7  5
\putrule from 10 4 to 12 4
\putrule from 12 4 to 12 6
\putrule from 10 6 to 12 6
\putrule from 10 4 to 10 6
\arrow <6pt> [.3,.67] from   11.09  4   to  11.11 4 
\put {$_{a}$}      at   11  3.7
\arrow <6pt> [.3,.67] from   12  5.09   to  12 5.11
\put {$_{b}$}      at   12.3  5
\arrow <6pt> [.3,.67] from   10.91  6   to  10.89 6 
\put {$_{\phantom{'}a'}$}      at   11  6.3
\arrow <6pt> [.3,.67] from   10  5.09   to  10 5.11
\put {$_{\phantom{'}b_{\beta + 1}}$}      at   9.5  5
\putrule from 6 0 to 8 0
\putrule from 8 0 to 8 2
\putrule from 6 2 to 8 2
\putrule from 6 0 to 6 2
\arrow <6pt> [.3,.67] from   7.09  0   to  7.11 0 
\put {$_{a}$}      at   7  -0.3
\arrow <6pt> [.3,.67] from   8  0.91   to  8 0.89
\put {$_{b_{\beta + 1}\phantom{'}}$}      at   8.6  1
\arrow <6pt> [.3,.67] from   6.91  2   to  6.89 2 
\put {$_{\phantom{'}a'}$}      at   7  2.3
\arrow <6pt> [.3,.67] from   6  0.91   to  6 0.89
\put {$_{b'}$}      at   5.7  1
\putrule from 10 0 to 12 0
\putrule from 12 0 to 12 2
\putrule from 10 2 to 12 2
\putrule from 10 0 to 10 2
\arrow <6pt> [.3,.67] from   11.09  0   to  11.11 0 
\put {$_{a}$}      at   11  -0.3
\arrow <6pt> [.3,.67] from   12  1.09   to  12 1.11
\put {$_{b}$}      at   12.3  1
\arrow <6pt> [.3,.67] from   10.91  2   to  10.89 2 
\put {$_{\phantom{'}a'}$}      at   11  2.3
\arrow <6pt> [.3,.67] from   10  0.91   to  10 0.89
\put {$_{\phantom{'}b_{\beta + 1}}$}      at   9.5  1
\putrule from 0 2 to 2 2
\putrule from 2 2 to 2 4
\putrule from 0 4 to 2 4
\putrule from 0 2 to 0 4
\arrow <6pt> [.3,.67] from   1.09  2   to  1.11 2 
\put {$_{a}$}      at   1  1.7
\arrow <6pt> [.3,.67] from   2  3.09   to  2 3.11
\put {$_{b}$}      at   2.3  3
\arrow <6pt> [.3,.67] from   0.91  4   to  0.89 4 
\put {$_{\phantom{'}a'}$}      at   1  4.3
\arrow <6pt> [.3,.67] from   0  2.91   to  0 2.89
\put {$_{b'}$}      at   -0.3  3
\arrow <6pt> [.3,.67] from   3  3.1   to  5 4.9
\arrow <6pt> [.3,.67] from   3  2.9   to  5 1.1
\endpicture
}
\hfil
\caption{The map $\phi_{\beta}$ associates
to the geometric square $[aba'b'] \in GS_{1,\beta}$ 
(represented on the left) 
the two geometric squares in $GS_{1,\beta + 1}$ represented on top right, 
and the two geometric squares in $GS_{1,\beta + 1}$ represented on bottom right, 
respectively.}
\end{figure}

\begin{Lemma} \label{Lemma1}
The map $\phi_{\beta}$ is well-defined.
\end{Lemma}
\begin{proof}
We have to show
\[
\phi_{\beta}([aba'b']) 
= \phi_{\beta}([a'b'ab]) 
= \phi_{\beta}([a^{-1}b'^{-1}a'^{-1}b^{-1}]) 
= \phi_{\beta}([a'^{-1}b^{-1}a^{-1}b'^{-1}]).
\]
Let $v_1 := [a b_{\beta + 1} a' b']$, 
$v_2 := [a b a' b_{\beta + 1}^{-1}]$,
$v_3 := [a b_{\beta + 1}^{-1} a' b']$
and
$v_4 := [a b a' b_{\beta + 1}]$,
such that we have $\{ v_1, v_2, v_3, v_4 \} \subset GS_{1,\beta + 1}$ and
$\phi_{\beta}([aba'b']) = \big\{ \{ v_1, v_2 \}, \{ v_3, v_4 \} \big\}$.

We check that
\begin{align}
\phi_{\beta}([a'b'ab]) &= \big\{ \{ [a' b_{\beta + 1} a b], [a' b' a b_{\beta + 1}^{-1}] \},
\{ [a' b_{\beta + 1}^{-1} a b], [a' b' a b_{\beta + 1}] \} \big\} \notag \\
&= \big\{ \{ v_4, v_3 \}, \{ v_2, v_1 \} \big\} = \phi_{\beta}([aba'b']), \notag
\end{align}
\begin{align}
\phi_{\beta}([a^{-1}b'^{-1}a'^{-1}b^{-1}]) 
= \big\{ &\{ [a^{-1} b_{\beta + 1} a'^{-1} b^{-1}], [a^{-1} b'^{-1} a'^{-1} b_{\beta + 1}^{-1}] \}, \notag \\
&\{ [a^{-1} b_{\beta + 1}^{-1} a'^{-1} b^{-1}], [a^{-1} b'^{-1} a'^{-1} b_{\beta + 1}] \} \big\} \notag \\
= \big\{ &\{ v_2, v_1 \}, \{ v_4, v_3 \} \big\} = \phi_{\beta}([aba'b']), \notag
\end{align}
\begin{align}
\phi_{\beta}([a'^{-1}b^{-1}a^{-1}b'^{-1}]) 
= \big\{ &\{ [a'^{-1} b_{\beta + 1} a^{-1} b'^{-1}], [a'^{-1} b^{-1} a^{-1} b_{\beta + 1}^{-1}] \}, \notag \\
&\{ [a'^{-1} b_{\beta + 1}^{-1} a^{-1} b'^{-1}], [a'^{-1} b^{-1} a^{-1} b_{\beta + 1}] \} \big\} \notag \\
= \big\{ &\{ v_3, v_4 \}, \{ v_1, v_2 \} \big\} = \phi_{\beta}([aba'b']). \notag
\end{align}
\end{proof} 
See Figure~\ref{fig2} for a visualization of the map $\phi_{\beta}$.

We now construct the set $\psi_{\beta}^{(2)}(R)$ consisting
of $2 \beta$ distinct $(1,\beta + 1)$-BM relations (as we will prove later).
Let
\[
\psi_{\beta}^{(2)}(R) := \bigcup_{i=1}^{\beta}
\left( \bigcup_{P \in \phi_{\beta}(r_i)} 
\big\{ P \cup (R \setminus \{ r_i \})\big\} \right).
\]
Note that if $r_i = [aba'b']$ then by definition of $\phi_{\beta}$
\[
\bigcup_{P \in \phi_{\beta}(r_i)} 
\big\{ P \cup (R \setminus \{ r_i \})\big\} =
\]
\[
\big\{ \{ [a b_{\beta + 1} a' b'], [a b a' b_{\beta + 1}^{-1}] \} \cup
(R \setminus \{r_i\}), \,
\{ [a b_{\beta + 1}^{-1} a' b'], [a b a' b_{\beta + 1}] \} \cup
(R \setminus \{r_i\})
\big\}.
\]

Finally, let 
\[
\psi_{\beta}(R) := \psi_{\beta}^{(1)}(R) \cup \psi_{\beta}^{(2)}(R).
\]
See Section~\ref{Appendix} for an explicit construction of the map $\psi_{\beta}$
in the case $\beta = 1$ and $\beta = 2$. 

\begin{Lemma} \label{Lemma2}
If $R \in R_{1,\beta}$, then the elements in $\psi_{\beta}(R)$ are 
$(1,\beta + 1)$-BM relations.
\end{Lemma}
\begin{proof}
The statement is clear for the three elements in $\psi_{\beta}^{(1)}(R)$ looking at their link.

To show it for the elements in $\psi_{\beta}^{(2)}(R)$,
first note that 
\[
[a b_{\beta + 1} a' b'] \ne [a b a' b_{\beta + 1}^{-1}]
\; (=[a^{-1} b_{\beta + 1} a'^{-1} b^{-1}])
\]
and 
\[
[a b_{\beta + 1}^{-1} a' b'] \ne [a b a' b_{\beta + 1}]
\; (=[a^{-1} b_{\beta + 1}^{-1} a'^{-1} b^{-1}]).
\]
Therefore each element in $\psi_{\beta}^{(2)}(R)$ consists of $\beta + 1$
geometric squares in $GS_{1,\beta + 1}$.
Let $R = \{r_1, \ldots , r_{\beta}\} \in R_{1,\beta}$, fix any 
$i \in \{ 1, \ldots, \beta \}$, and suppose that
$r_i = [aba'b']$.
Since $Lk(R) = K_{2,2\beta}$,
we have
\[
Lk(\{r_1, \ldots , r_{\beta}, 
[a b_{\beta + 1} a' b_{\beta + 1}^{-1}] \}) = K_{2,2\beta + 2},
\]
independently of $a, a' \in \{a_1, a_1^{-1}\} = A_1^{\pm 1}$.
Since $r_i = [aba'b']$, we can write this as
\[
K_{2,2\beta + 2} = 
Lk(\{ [a b a' b'], [a b_{\beta + 1} a' b_{\beta + 1}^{-1}] \} \cup
(R \setminus \{r_i\})),
\]
which can directly be seen to be equal to
\[
Lk(\{ [a b_{\beta + 1} a' b'], [a b a' b_{\beta + 1}^{-1}] \} \cup
(R \setminus \{r_i\})),
\]
since the edges in
$Lk(\{ [a b a' b'], [a b_{\beta + 1} a' b_{\beta + 1}^{-1}] \})$
are 
\begin{align}
&\{ a^{-1}, b  \}, \{ a', b^{-1}  \}, 
\{ a'^{-1}, b'  \}, \{ a, b'^{-1} \}, \notag \\
&\{ a^{-1}, b_{\beta + 1}  \}, \{ a', b_{\beta + 1}^{-1}  \}, 
\{ a'^{-1}, b_{\beta + 1}^{-1}  \}, \{ a, b_{\beta + 1}  \}, \notag
\end{align}
which are also the edges in
$Lk(\{ [a b_{\beta + 1} a' b'], [a b a' b_{\beta + 1}^{-1}] \})$.
In fact, we have performed a link preserving surgery as described more generally in
\cite[Section~6.2.2]{BMII}.

Similarly (interchanging $b_{\beta + 1}$ and $b_{\beta + 1}^{-1}$) one proves that
\[
Lk(\{ [a b_{\beta + 1}^{-1} a' b'], [a b a' b_{\beta + 1}] \} \cup
(R \setminus \{r_i\})) = K_{2,2\beta + 2}.
\]
\end{proof}

\begin{Lemma} \label{Lemma3}
If $R \in R_{1,\beta}$, then $|\psi_{\beta}(R)| = 3 + 2 \beta$.
\end{Lemma}
\begin{proof}
Clearly $|\psi_{\beta}^{(1)}(R)| = 3$.

The label $b_{\beta + 1}$ or $b_{\beta + 1}^{-1}$ appears in exactly 
one geometric square of each element in $\psi_{\beta}^{(1)}(R)$,
but in exactly two geometric squares of each element in $\psi_{\beta}^{(2)}(R)$,
hence we conclude 
\[
\psi_{\beta}^{(1)}(R) \cap \psi_{\beta}^{(2)}(R) = \emptyset.
\]

Let $R = \{r_1, \ldots , r_{\beta}\}$.
Fix any $i \in \{1, \ldots, \beta \}$
and suppose that $r_i = [aba'b']$.
The geometric square $r_i$ only misses in the two elements
\[
\{ [a b_{\beta + 1} a' b'], [a b a' b_{\beta + 1}^{-1}] \} \cup (R \setminus \{r_i\})
\]
and
\[
\{ [a b_{\beta + 1}^{-1} a' b'], [a b a' b_{\beta + 1}] \} \cup (R \setminus \{r_i\})
\]
of $\psi_{\beta}^{(2)}(R)$.
Suppose that they are equal. Then
\[
\{ [a b_{\beta + 1} a' b'], [a b a' b_{\beta + 1}^{-1}] \}
= \{ [a b_{\beta + 1}^{-1} a' b'], [a b a' b_{\beta + 1}] \}.
\]
It follows that 
\[
[a b_{\beta + 1} a' b'] = [a b a' b_{\beta + 1}] 
\; (=[a' b_{\beta + 1} a b]),
\]
since $[a b_{\beta + 1} a' b'] \ne [a b_{\beta + 1}^{-1} a' b']$,
but then $a=a'$ and $b=b'$.
This is impossible, since $[abab] \notin R \in R_{1,\beta}$.
This shows that the $2 \beta$ elements in $\psi_{\beta}^{(2)}(R)$ are distinct, and
we get 
\[
|\psi_{\beta}(R)| = |\psi_{\beta}^{(1)}(R) \cup \psi_{\beta}^{(2)}(R)| = 
|\psi_{\beta}^{(1)}(R)| + |\psi_{\beta}^{(2)}(R)| - |\psi_{\beta}^{(1)}(R) \cap \psi_{\beta}^{(2)}(R)|
= 3 + 2 \beta.
\]
\end{proof}

\begin{Lemma} \label{Lemma4}
If $R, T \in R_{1,\beta}$ and $R \ne T$, then $\psi_{\beta}(R) \cap \psi_{\beta}(T) = \emptyset$.
\end{Lemma}
\begin{proof}
Let $R = \{r_1, \ldots , r_{\beta}\}$ and 
$T = \{t_1, \ldots , t_{\beta}\}$.
We suppose without loss of generality that
$r_1 = [aba'b'] \notin T$.
Then $r_1$ appears in no element of $\psi_{\beta}(T)$,
but appears in each element of $\psi_{\beta}(R)$
except in 
\[
U_1 := \{ [a b_{\beta + 1} a' b'], [a b a' b_{\beta + 1}^{-1}] \} \cup
(R \setminus \{r_1\}) \in \psi_{\beta}^{(2)}(R)
\]
and
\[
V_1 := \{ [a b_{\beta + 1}^{-1} a' b'], [a b a' b_{\beta + 1}] \} \cup
(R \setminus \{r_1\}) \in \psi_{\beta}^{(2)}(R).
\]
We want to show by contradiction that $U_1, V_1 \notin \psi_{\beta}(T)$.
It is clear that $U_1, V_1 \notin \psi_{\beta}^{(1)}(T)$.  
Fix any $i \in \{1, \dots, \beta \}$ and let 
$t_i = [\breve{a} \breve{b} \hat{a} \hat{b}]$,
where $\breve{a}, \hat{a} \in \{ a_1, a_1^{-1} \}$
and $\breve{b}, \hat{b} \in B_{\beta}^{\pm 1}$.
We suppose that $U_1 \in \psi_{\beta}^{(2)}(T)$ or $V_1 \in \psi_{\beta}^{(2)}(T)$
and have therefore to consider four cases:

Case~1:
Suppose that 
\[
U_1
= \{ [\breve{a} b_{\beta + 1} \hat{a} \hat{b}], 
[\breve{a} \breve{b} \hat{a} b_{\beta + 1}^{-1}] \} \cup (T \setminus \{t_i\}).
\]
Then $R \setminus \{r_1\} = T \setminus \{t_i\}$
and 
\[
\{ [a b_{\beta + 1} a' b'], [a b a' b_{\beta + 1}^{-1}] \}
= \{ [\breve{a} b_{\beta + 1} \hat{a} \hat{b}], 
[\breve{a} \breve{b} \hat{a} b_{\beta + 1}^{-1}] \}.
\]

Case~1.1:
If $[a b_{\beta + 1} a' b'] = 
[\breve{a} b_{\beta + 1} \hat{a} \hat{b}]$
and $[a b a' b_{\beta + 1}^{-1}] =
[\breve{a} \breve{b} \hat{a} b_{\beta + 1}^{-1}]$,
then $a = \breve{a}$, $b = \breve{b}$, $a' = \hat{a}$
and $b' = \hat{b}$.
This implies $t_i = [aba'b'] = r_1$, hence $R=T$,
a contradiction. 

Case~1.2:
If 
\[
[a b_{\beta + 1} a' b'] =
[\breve{a} \breve{b} \hat{a} b_{\beta + 1}^{-1}]
\; (=[\breve{a}^{-1} b_{\beta + 1} \hat{a}^{-1} \breve{b}^{-1}])
\]
and 
\[
[a b a' b_{\beta + 1}^{-1}] =
[\breve{a} b_{\beta + 1} \hat{a} \hat{b}]
\; (=[\breve{a}^{-1} \hat{b}^{-1} \hat{a}^{-1} b_{\beta + 1}^{-1}]),
\]
then $a = \breve{a}^{-1}$, $b = \hat{b}^{-1}$, $a' = \hat{a}^{-1}$
and $b' = \breve{b}^{-1}$.
This implies 
\[
t_i = [a^{-1}b'^{-1}a'^{-1}b^{-1}] = [aba'b'] = r_1
\]
and again the contradiction $R=T$.

The three remaining cases

Case 2: 
$U_1
= \{ [\breve{a} b_{\beta + 1}^{-1} \hat{a} \hat{b}], 
[\breve{a} \breve{b} \hat{a} b_{\beta + 1}] \} \cup (T \setminus \{t_i\})$

Case 3: 
$V_1
= \{ [\breve{a} b_{\beta + 1} \hat{a} \hat{b}], 
[\breve{a} \breve{b} \hat{a} b_{\beta + 1}^{-1}] \} \cup (T \setminus \{t_i\})$

Case 4:
$V_1
= \{ [\breve{a} b_{\beta + 1}^{-1} \hat{a} \hat{b}], 
[\breve{a} \breve{b} \hat{a} b_{\beta + 1}] \} \cup (T \setminus \{t_i\})$

\noindent can be treated similarly. In fact we can reduce them to Case~1 as follows:

In Case~2, since 
\[
\{ [\breve{a} b_{\beta + 1}^{-1} \hat{a} \hat{b}], 
[\breve{a} \breve{b} \hat{a} b_{\beta + 1}] \}
=
\{ [\hat{a}^{-1} b_{\beta + 1} \breve{a}^{-1} \hat{b}^{-1}], 
[\hat{a}^{-1} \breve{b}^{-1} \breve{a}^{-1} b_{\beta + 1}^{-1}] \},
\]
we can substitute $\breve{a} \breve{b} \hat{a} \hat{b}$ by 
$\hat{a}^{-1} \breve{b}^{-1} \breve{a}^{-1} \hat{b}^{-1}$
and are in Case~1.

In Case~3 and Case~4, since
\[
(V_1 \cup \{r_1\}) \setminus R = \{ [a b_{\beta + 1}^{-1} a' b'], [a b a' b_{\beta + 1}] \} 
= \{ [a'^{-1} b_{\beta + 1} a^{-1} b'^{-1}], [a'^{-1} b^{-1} a^{-1} b_{\beta + 1}^{-1}] \},
\]
we can substitute $aba'b'$ by $a'^{-1} b^{-1} a^{-1} b'^{-1}$ and are
in Case~1 and Case~2, respectively.

Thus, we have shown that the only two elements $U_1, V_1$ of $\psi_{\beta}(R)$ 
in which $r_1$ does not appear,
are no elements of $\psi_{\beta}(T)$,
and therefore $\psi_{\beta}(R) \cap \psi_{\beta}(T) = \emptyset$.
\end{proof}

\begin{Lemma} \label{Lemma5}
Let $U \in R_{1,\beta + 1}$. 
Then $U \in \psi_{\beta}(R)$ for some $R \in R_{1,\beta}$.
\end{Lemma}
\begin{proof}
Let $U = \{ u_1, \ldots, u_{\beta + 1}\} \in R_{1,\beta + 1}$.
By the link condition, the label $b_{\beta + 1}$ or $b_{\beta + 1}^{-1}$
appears either in exactly one or in exactly two elements (geometric squares) of $U$.

Case 1: Suppose that the label $b_{\beta + 1}$ or $b_{\beta + 1}^{-1}$
appears in exactly one element of $U$, say in $u_{\beta + 1}$.
Then either
\[
u_{\beta + 1} = [a_1 b_{\beta + 1} a_1^{-1} b_{\beta + 1}^{-1}]
\]
or
\[
u_{\beta + 1} = [a_1 b_{\beta + 1} a_1 b_{\beta + 1}^{-1}]
\]
or
\[
u_{\beta + 1} = [a_1 b_{\beta + 1} a_1^{-1} b_{\beta + 1}].
\]
Let $R := \{ u_1, \ldots, u_{\beta}\} = U \setminus \{ u_{\beta + 1} \}$. 
Note that $R \in R_{1,\beta}$, since $Lk(U) = K_{2,2\beta + 2}$
and $u_{\beta + 1}$ contributes to $Lk(U)$ the four edges 
\[
\{a_1^{-1}, b_{\beta + 1} \}, \{a_1^{-1}, b_{\beta + 1}^{-1} \}, 
\{a_1, b_{\beta + 1}^{-1} \}, \{a_1, b_{\beta + 1} \},
\]
independently of the three possibilities for $u_{\beta + 1}$.
By definition of $\psi_{\beta}^{(1)}$ we have 
$U \in \psi_{\beta}^{(1)}(R) \subset \psi_{\beta}(R)$.

Case 2:  Suppose that the label $b_{\beta + 1}$ or $b_{\beta + 1}^{-1}$
appears in $u_{\beta}$ and $u_{\beta + 1}$, but in no other element of $U$.
It follows that $u_{\beta} = [ab_{\beta + 1} a' b']$ for some
$a, a' \in \{a_1, a_1^{-1}\}$ and $b' \in B_{\beta}^{\pm 1}$.
In particular $b' \ne b_{\beta + 1}$ and $b' \ne b_{\beta + 1}^{-1}$,
otherwise we would be in Case~1.
Looking at the link of 
$\{ u_1, \ldots, u_{\beta}\} = U \setminus \{u_{\beta + 1}\}$,
we see that the two edges $\{ a^{-1}, b_{\beta + 1}\}$ and 
$\{a', b_{\beta + 1}^{-1}\}$ in this link
(and two other edges not involving the label $b_{\beta + 1}$ or $b_{\beta + 1}^{-1}$)
are contributed by $u_{\beta}$. The edges contributed by 
$\{ u_1, \ldots, u_{\beta - 1}\}$ do not involve $b_{\beta + 1}$ or $b_{\beta + 1}^{-1}$.
Therefore, the two edges $\{a, b_{\beta + 1}\}$ and 
$\{a'^{-1}, b_{\beta + 1}^{-1}\}$ 
(and two other edges not involving the label $b_{\beta + 1}$ or $b_{\beta + 1}^{-1}$) 
are missing to get the complete bipartite graph $K_{2,2\beta + 2} = Lk(U)$.
Hence $u_{\beta + 1} = [aba' b_{\beta + 1}^{-1}]$ for some
$b \in B_{\beta}^{\pm 1}$.
Let $R := \{ u_1, \ldots, u_{\beta - 1}, [aba'b']\}$.
Then $R \in R_{1,\beta}$ (i.e.\ $Lk(R) = K_{2,2\beta}$), since 
\begin{align}
K_{2,2\beta + 2} = Lk(U) &= Lk(\{ u_1, \ldots, u_{\beta - 1},
[ab_{\beta + 1} a' b'], [aba' b_{\beta + 1}^{-1}]\}) \notag \\
&= Lk(\{ u_1, \ldots, u_{\beta - 1}, 
[aba'b'], [ab_{\beta + 1} a' b_{\beta + 1}^{-1}]\}) \notag \\
&= Lk(R \cup \{ [ab_{\beta + 1} a' b_{\beta + 1}^{-1}] \}). \notag
\end{align}
By construction of $R$ and the definition of $\psi_{\beta}^{(2)}$, 
we have $U \in \psi_{\beta}^{(2)}(R) \subset \psi_{\beta}(R)$.
\end{proof}

\begin{Corollary} \label{CorUnion}
For $\beta \in \mathbb{N}$ we have
\[
\bigcup_{R \in R_{1,\beta}} \psi_{\beta}(R) = R_{1,\beta + 1},
\]
in particular the set $R_{1,\beta + 1}$ can be explicitly constructed from
$R_{1,\beta}$ using $\psi_{\beta}$.
\end{Corollary}
\begin{proof}
Lemma~\ref{Lemma2} shows that
\[
\bigcup_{R \in R_{1,\beta}} \psi_{\beta}(R) \subseteq R_{1,\beta + 1}.
\]
Moreover, we have
\[
\bigcup_{R \in R_{1,\beta}} \psi_{\beta}(R) \supseteq R_{1,\beta + 1}
\]
by Lemma~\ref{Lemma5}.
\end{proof}
Note that the union in Corollary~\ref{CorUnion} is a \emph{disjoint} union by Lemma~\ref{Lemma4}.
Now, we are able to prove Kimberley's conjecture on
the number of $(1,\beta )$-BM relations. 
\begin{Theorem} \label{Theorem193}
(\cite[Conjecture~193]{JSK})
For every positive integer $\beta$, the number of 
$(1,\beta )$-BM relations is
\[
|R_{1,\beta}| = \prod_{i=1}^{\beta} (2i+1).
\]
\end{Theorem}
\begin{proof}
By Lemma~\ref{Lemma3} and Lemma~\ref{Lemma5}
\[
|R_{1,\beta + 1}| \leq (3+2\beta)|R_{1,\beta}|.
\]
By Lemma~\ref{Lemma3} and Lemma~\ref{Lemma4}
\[
|R_{1,\beta + 1}| \geq (3+2\beta)|R_{1,\beta}|,
\]
hence 
\[
|R_{1,\beta + 1}| = (3+2\beta)|R_{1,\beta}|.
\]
The proof of the theorem is now by induction on $\beta$.
If $\beta =1$, then 
\[
R_{1,1} = \big\{ \{[a_1 b_1 a_1^{-1} b_1^{-1}]\}, 
\{[a_1 b_1 a_1 b_1^{-1}]\},
\{[a_1 b_1 a_1^{-1} b_1]\} \big\}
\]
and $|R_{1,1}| = 3$. Assume that the statement of the theorem
holds for $\beta$. Then
\[
|R_{1,\beta + 1}| = (3+2\beta)|R_{1,\beta}| = 
(2(\beta + 1) + 1) \prod_{i=1}^{\beta} (2i+1) = \prod_{i=1}^{\beta + 1} (2i+1).
\]
\end{proof}

\section{Classification of $(2,2)$-BM groups} \label{SectionClass}
Let $\Gamma_{4}$, $\Gamma_{30}$, $\Gamma_{5}$, $\Gamma_{10}$
be the $(2,2)$-BM groups 
\[
\Gamma_{4} = \langle a,b,c,d \mid acac^{-1}, \, adad^{-1}, \, bcbd, \, bc^{-1}bd^{-1} \rangle,
\]
\[
\Gamma_{30} = \langle a,b,c,d \mid acad, \, ac^{-1}ad^{-1}, \, bcbd, \, bc^{-1}bd^{-1} \rangle,
\]
\[
\Gamma_{5} = \langle a,b,c,d \mid acac^{-1}, \, adad^{-1}, \, bcb^{-1}c, \, bdb^{-1}d \rangle,
\]
\[
\Gamma_{10} = \langle a,b,c,d \mid acac^{-1}, \, ada^{-1}d, \, bcbc^{-1}, \, bdb^{-1}d^{-1} \rangle.
\]
(To simplify the notation, we use here the letters $a, b, c, d$ 
instead of $a_1, a_2, b_1, b_2$.)

We will prove that $\Gamma_4$ is isomorphic to $\Gamma_{30}$, and that
$\Gamma_5$ is isomorphic to $\Gamma_{10}$.
To find these isomorphisms we have written a program with GAP (\cite{GAP})
using the normal form program developed in \cite[Chapter~B.6]{Rattaggi}
and the knowledge of the orders of elements in the abelianizations
of the four groups.

\begin{Proposition} \label{Prop1}
The groups $\Gamma_4$ and $\Gamma_{30}$ are isomorphic.
\end{Proposition}
\begin{proof}
Let $\eta : \Gamma_4 \to \Gamma_{30}$ be the homomorphism given by 
$\eta(a)=ab$, $\eta(b)=a$, $\eta(c)=ac$ and $\eta(d)=da^{-1}$. 
It is a homomorphism since 
\[
\eta(acac^{-1}) = abacabc^{-1}a^{-1} = abaa^{-1}d^{-1}bad =
abd^{-1}bad = acad = 1,
\] 
\[
\eta(adad^{-1}) = abda^{-1}abad^{-1} = abdbad^{-1} =
abb^{-1}c^{-1}ad^{-1} = ac^{-1}ad^{-1} = 1,
\]
\[
\eta(bcbd) = aacada^{-1} = aaa^{-1}d^{-1}da^{-1} = 1,
\]
\[
\eta(bc^{-1}bd^{-1}) = ac^{-1}a^{-1}aad^{-1} = ac^{-1}ad^{-1} = 1,
\]
using the four defining relations of $\Gamma_{30}$.
 
$\eta$ is surjective: $a = \eta(b)$, $b = \eta(b^{-1}a)$, $c =
\eta(b^{-1}c)$, $d = \eta(db)$. 

Let $\theta : \Gamma_{30} \to \Gamma_{4}$ be the homomorphism given by 
$\theta(a)=b$, $\theta(b)=b^{-1}a$, $\theta(c)=b^{-1}c$ and $\theta(d)=db$. 
It is a homomorphism since
\[
\theta(acad) = bb^{-1}cbdb = cbdb = 1,
\]
\[
\theta(ac^{-1}ad^{-1}) = bc^{-1}bbb^{-1}d^{-1} = bc^{-1}bd^{-1} = 1,
\]
\[
\theta(bcbd) = b^{-1}ab^{-1}cb^{-1}adb = b^{-1}ab^{-1}bd^{-1}da^{-1}b = 1,
\]
\[
\theta(bc^{-1}bd^{-1}) = b^{-1}ac^{-1}bb^{-1}ab^{-1}d^{-1} =
b^{-1}c^{-1}a^{-1}ab^{-1}d^{-1} = b^{-1}c^{-1}b^{-1}d^{-1} = 1,
\]
using the four defining relations of $\Gamma_{4}$.

The composition $\theta \circ \eta$ is the identity on $\Gamma_4$, since
\[
\theta(\eta(a)) = \theta(ab) = bb^{-1}a = a, 
\]
\[
\theta(\eta(b)) = \theta(a) = b,
\]
\[
\theta(\eta(c)) = \theta(ac) = bb^{-1}c = c, 
\]
\[
\theta(\eta(d)) = \theta(da^{-1}) = dbb^{-1} = d, 
\]
hence $\eta$ is injective and an isomorphism.  
\end{proof}

\begin{Proposition} \label{Prop2}
The groups $\Gamma_5$ and $\Gamma_{10}$ are isomorphic.
\end{Proposition}
\begin{proof}
As in the proof of Proposition~\ref{Prop1},
it is easy to show that $\varphi : \Gamma_{5} \to \Gamma_{10}$ 
defined by
$\varphi(a)=d$, $\varphi(b)=ac$, $\varphi(c)=a$, $\varphi(d)=ab$
is an isomorphism.
\end{proof}

\begin{Corollary}
There are exactly $41$ $(2,2)$-BM groups up to isomorphism.
\end{Corollary}
\begin{proof}
By \cite[Proposition~222]{JSK} there are at least $41$ isomorphism
classes of $(2,2)$-BM groups.
By \cite[Proposition~231]{JSK} there are at most $43$ isomorphism
classes of $(2,2)$-BM groups
(including the isomorphism classes of $\Gamma_{4}$, $\Gamma_{30}$,
$\Gamma_{5}$ and $\Gamma_{10}$).
Now use Proposition~\ref{Prop1} and Proposition~\ref{Prop2}
to reduce the number of isomorphism classes from $43$ to $41$.
\end{proof}

\section{Appendix: Illustration of $\psi_{\beta}$ for $\beta = 1$ and $\beta = 2$} \label{Appendix}
In this appendix we first use the map $\psi_{1}$ to determine 
\[
\bigcup_{R \in R_{1,1}} \psi_{1}(R) = R_{1,2}.
\]
Recall that
\[
R_{1,1} = \big\{ \{[a_1 b_1 a_1^{-1} b_1^{-1}]\}, 
\{[a_1 b_1 a_1 b_1^{-1}]\},
\{[a_1 b_1 a_1^{-1} b_1]\} \big\}.
\]
By definition of $\psi_{1}^{(1)}$ and $\psi_{1}^{(2)}$, we have
\begin{align}
\psi_{1}^{(1)}(\{[a_1 b_1 a_1^{-1} b_1^{-1}]\}) = \big\{ 
&\{ [a_1 b_1 a_1^{-1} b_1^{-1}], [a_1 b_2 a_1^{-1} b_2^{-1}] \}, \notag \\
&\{ [a_1 b_1 a_1^{-1} b_1^{-1}], [a_1 b_2 a_1 b_2^{-1}] \}, \notag \\
&\{ [a_1 b_1 a_1^{-1} b_1^{-1}], [a_1 b_2 a_1^{-1} b_2] \} \big\}. \notag \\
\psi_{1}^{(2)}(\{[a_1 b_1 a_1^{-1} b_1^{-1}]\}) = \big\{
&\{ [a_1 b_2 a_1^{-1} b_1^{-1}], [a_1 b_1 a_1^{-1} b_2^{-1}] \}, \notag \\
&\{ [a_1 b_2^{-1} a_1^{-1} b_1^{-1}], [a_1 b_1 a_1^{-1} b_2] \} \big\}. \notag 
\end{align}
\begin{align}
\psi_{1}^{(1)}(\{[a_1 b_1 a_1 b_1^{-1}]\}) = \big\{ 
&\{ [a_1 b_1 a_1 b_1^{-1}], [a_1 b_2 a_1^{-1} b_2^{-1}] \}, \notag \\
&\{ [a_1 b_1 a_1 b_1^{-1}], [a_1 b_2 a_1 b_2^{-1}] \}, \notag \\
&\{ [a_1 b_1 a_1 b_1^{-1}], [a_1 b_2 a_1^{-1} b_2] \} \big\}. \notag \\
\psi_{1}^{(2)}(\{[a_1 b_1 a_1 b_1^{-1}]\}) = \big\{
&\{ [a_1 b_2 a_1 b_1^{-1}], [a_1 b_1 a_1 b_2^{-1}] \}, \notag \\
&\{ [a_1 b_2^{-1} a_1 b_1^{-1}], [a_1 b_1 a_1 b_2] \} \big\}. \notag 
\end{align}
\begin{align}
\psi_{1}^{(1)}(\{[a_1 b_1 a_1^{-1} b_1]\}) = \big\{ 
&\{ [a_1 b_1 a_1^{-1} b_1], [a_1 b_2 a_1^{-1} b_2^{-1}] \}, \notag \\
&\{ [a_1 b_1 a_1^{-1} b_1], [a_1 b_2 a_1 b_2^{-1}] \}, \notag \\
&\{ [a_1 b_1 a_1^{-1} b_1], [a_1 b_2 a_1^{-1} b_2] \} \big\}. \notag \\
\psi_{1}^{(2)}(\{[a_1 b_1 a_1^{-1} b_1]\}) = \big\{
&\{ [a_1 b_2 a_1^{-1} b_1], [a_1 b_1 a_1^{-1} b_2^{-1}] \}, \notag \\
&\{ [a_1 b_2^{-1} a_1^{-1} b_1], [a_1 b_1 a_1^{-1} b_2] \} \big\}. \notag 
\end{align} 
Taking the union of these six sets, we therefore obtain
\begin{align}
R_{1,2} = \big\{ 
&\{ [a_1 b_1 a_1^{-1} b_1^{-1}], [a_1 b_2 a_1^{-1} b_2^{-1}] \}, \,
\{ [a_1 b_1 a_1^{-1} b_1^{-1}], [a_1 b_2 a_1 b_2^{-1}] \}, \notag \\
&\{ [a_1 b_1 a_1^{-1} b_1^{-1}], [a_1 b_2 a_1^{-1} b_2] \},  \,
\{ [a_1 b_2 a_1^{-1} b_1^{-1}], [a_1 b_1 a_1^{-1} b_2^{-1}] \}, \notag \\
&\{ [a_1 b_2^{-1} a_1^{-1} b_1^{-1}], [a_1 b_1 a_1^{-1} b_2] \}, \,
\{ [a_1 b_1 a_1 b_1^{-1}], [a_1 b_2 a_1^{-1} b_2^{-1}] \}, \notag \\
&\{ [a_1 b_1 a_1 b_1^{-1}], [a_1 b_2 a_1 b_2^{-1}] \}, \,
\{ [a_1 b_1 a_1 b_1^{-1}], [a_1 b_2 a_1^{-1} b_2] \},  \notag \\
&\{ [a_1 b_2 a_1 b_1^{-1}], [a_1 b_1 a_1 b_2^{-1}] \}, \,
\{ [a_1 b_2^{-1} a_1 b_1^{-1}], [a_1 b_1 a_1 b_2] \}, \notag \\
&\{ [a_1 b_1 a_1^{-1} b_1], [a_1 b_2 a_1^{-1} b_2^{-1}] \}, \,
\{ [a_1 b_1 a_1^{-1} b_1], [a_1 b_2 a_1 b_2^{-1}] \}, \notag \\
&\{ [a_1 b_1 a_1^{-1} b_1], [a_1 b_2 a_1^{-1} b_2] \}, \,
\{ [a_1 b_2 a_1^{-1} b_1], [a_1 b_1 a_1^{-1} b_2^{-1}] \}, \notag \\
&\{ [a_1 b_2^{-1} a_1^{-1} b_1], [a_1 b_1 a_1^{-1} b_2] \} \notag
\big\} 
\end{align}
and $|R_{1,2}| = 15$. 
These $15$ $(1,2)$-BM relations are also listed in \cite[Table~7]{JSK}.

To illustrate what happens in the case $\beta = 2$, we take for example 
\[
R := \{ [a_1 b_1 a_1^{-1} b_1^{-1}], [a_1 b_2 a_1^{-1} b_2^{-1}] \} \in R_{1,2},
\]
and get seven $(1,3)$-BM relations
\begin{align}
\psi_{2}^{(1)}(R) = \big\{ 
&\{ [a_1 b_1 a_1^{-1} b_1^{-1}], [a_1 b_2 a_1^{-1} b_2^{-1}], [a_1 b_3 a_1^{-1} b_3^{-1}] \}, \notag \\
&\{ [a_1 b_1 a_1^{-1} b_1^{-1}], [a_1 b_2 a_1^{-1} b_2^{-1}], [a_1 b_3 a_1 b_3^{-1}] \}, \notag \\
&\{ [a_1 b_1 a_1^{-1} b_1^{-1}], [a_1 b_2 a_1^{-1} b_2^{-1}], [a_1 b_3 a_1^{-1} b_3] \} \big\}. \notag \\
\psi_{2}^{(2)}(R) = \big\{
&\{ [a_1 b_3 a_1^{-1} b_1^{-1}], [a_1 b_1 a_1^{-1} b_3^{-1}], 
[a_1 b_2 a_1^{-1} b_2^{-1}] \}, \notag \\
&\{ [a_1 b_3^{-1} a_1^{-1} b_1^{-1}], [a_1 b_1 a_1^{-1} b_3], 
[a_1 b_2 a_1^{-1} b_2^{-1}] \}, \notag \\
&\{ [a_1 b_3 a_1^{-1} b_2^{-1}], [a_1 b_2 a_1^{-1} b_3^{-1}], 
[a_1 b_1 a_1^{-1} b_1^{-1}] \}, \notag \\
&\{ [a_1 b_3^{-1} a_1^{-1} b_2^{-1}], [a_1 b_2 a_1^{-1} b_3], 
[a_1 b_1 a_1^{-1} b_1^{-1}] \} 
\big\}. \notag 
\end{align}


\begin{thebibliography}{9}
\bibitem{BMII}
  Burger, Marc; Mozes, Shahar,
  \emph{Lattices in product of trees},
  Inst. Hautes \'Etudes Sci. Publ. Math. No. \textbf{92}(2000),
  151--194(2001).
\bibitem{GAP}
  The GAP group, \textsf{GAP} --- Groups, Algorithms, and Programming, Version 4.4; 2004,
  \texttt{http://www.gap-system.org}
\bibitem{JSK}
  Kimberley, Jason S., 
  \emph{Classifying Burger-Mozes groups and the algebras generated
  from their actions},
  forthcoming Ph.D.\ thesis, University of Newcastle, Australia, 2006.
\bibitem{KR}
  Kimberley, Jason S.; Robertson, Guyan,
  \emph{Groups acting on products of trees, tiling systems and analytic K-theory},
  New York J. Math. \textbf{8}(2002), 111--131 (electronic).
\bibitem{Rattaggi}
  Rattaggi, Diego,
  \emph{Computations in groups acting on a product of trees: 
  normal subgroup structures and quaternion lattices},
  Ph.D.~thesis, ETH Z\"urich, 2004.
\bibitem{Rattaggi2} 
  Rattaggi, Diego,
  \emph{A finitely presented torsion-free simple group},
  to appear in J. Group Theory, 
  also available at arXiv:math.GR/0411546.
\bibitem{Wise}
  Wise, Daniel T.,
  \emph{Non-positively curved squared complexes, aperiodic tilings,
  and non-residually finite groups},
  Ph.D.~thesis, Princeton University, 1996.
\end{thebibliography}
\end{document}